\documentclass[12pt,reqno,oneside]{amsart}
\usepackage{amsmath,amsthm,amsfonts,amssymb}
\usepackage[mathscr]{eucal}
\usepackage{indentfirst}
\usepackage{url}

\theoremstyle{plain}
\newtheorem{teo}{Theorem}[section]

\newtheorem{lem}[teo]{Lemma}

\theoremstyle{definition}

\newtheorem{exa}[teo]{Example}

\numberwithin{equation}{section}

\def\bbR{{\mathbb R}}
\def\bbP{{\mathbb P}}
\def\bbZ{{\mathbb Z}}
\def\bbN{{\mathbb N}}
\def\bbE{{\mathbf E}}
\def\bbPr{{\mathbf P}}
\def\bbT{{\mathbb T}}

\def\qed{\hfill $\square$}

\def\cpp{\textit{Cone Percolation Process}}
\def\hcpp{\textit{Heterogeneous Cone Percolation Process}}

\begin{document}

\baselineskip=26pt

\title{The cone percolation on $\bbT_d$}

\author{Valdivino~V.~Junior}
\author{F\'abio~P.~Machado}
\author{Mauricio~Zuluaga}

\address[F\'abio~P.~Machado]
{Institute of Mathematics and Statistics
\\ University of S\~ao Paulo \\ Rua do Mat\~ao 1010, CEP
05508-090, S\~ao Paulo, SP, Brazil.}

\address[Valdivino~V.~Junior]
{Federal University of Goias
\\ Campus Samambaia, CEP 74001-970, Goi\^ania, GO, Brazil.}

\address[Mauricio~Zuluaga]
{Department of Statistics, Federal University of Pernambuco
\\ Cidade Universit\'aria, CEP 50740-540, Recife, PE, Brazil.}

\noindent
\email{vvjunior@mat.ufg.br, fmachado@ and zuluaga@ime.usp.br}

\thanks{Research supported by CNPq (306927/2007-1) and FAPESP (2010/50884-4 and 2009/18253-7).}



\date{\today}

\begin{abstract}
We study a rumour model from a percolation theory and branching process 
point of view. The existence of a giant component is related to the event
where the rumour spreads out trough an infinite number of individuals.
We present sharp lower and upper bounds for the
probability of that event, according to the distribution of the random variables 
that defines the radius of influence of each individual.
\end{abstract}

\maketitle

\section{Introduction}
\label{S: Introduction}

We study long range dependent oriented percolation processes on a tree through its 
most basic propriety: the existence of a giant connected component. 
The starting point for approaches to rigorous percolation theory beyond the 
nearest neighbor independent setup on $\bbZ^d$ is due to several authors around the nineties. 
Grimmett and Newman~\cite{GrimmettNewman} in 1990 study percolation
on $\bbT_d \times \bbZ$, Burton and Meester~\cite{BM} in 1993 study phase transition for a
long range independent percolation model on a stationary point process in $\bbR^d$, 
Lyons~\cite{Lyons} put out the first version of his book \textit{Probability
on Trees} in 1994 while Benjamini and Schram~\cite{BS} in 1996 have they \textit{Percolation beyond $\bbZ^d$} published,
just to name a few. Lebensztayn and Rodriguez~\cite{LR} in 2008, propose a model on general graphs 
named \textit{disk percolation} where a reaction chain starting from the origin of the graph, based on independent
copies of a geometric random variables with parameter $q \in [0,1],$ defines the existence or not of 
a giant component. They obtain a sufficient condition for the existence of phase transition based on $q$, which means 
the existence of a non-empty subcritical (no giant components) and supercritical (giant components with positive
probability) phases. They associate their model to a rumour or an epidemic process. In this paper, instead of working
in a general family of graphs we focus on homogeneous trees and instead of fixing the random variable which defines
the radius of infection or the radius of influence of each vertex to be geometric, we consider general random variables. So, as
a result, instead of having a phase transition phenomena dependending on a point in a parametric space, we have 
that phenomena depending on the family of general positive random variables. 


We consider a process which allows us to associate the activation dynamic 
on the set of vertices
to a discrete rumour process. In\-di\-vi\-duals become spreaders as soon as 
they heard about the rumour. Next time, they propagate the rumour within their 
\textit{radius of influence} and immediately become stiflers. Our main interest is to 
establish whether the process has positive probability of involving an 
infinite set of individual. Besides, we present sharp lower and upper bounds for the
probability of that event, according to the general distribution of the random variables 
that defines the \textit{radius of influence} of each individual. We say that 
the process \textit{survives} if the amount of vertices involved is infinite. 
Otherwise we say the process \textit{dies out}. 

Consider $\bbT_d$ the homogeneous tree such that each vertex has 
$d+1$ neighbours, $d \ge 2$. Let ${\mathcal V}(\bbT_d)$ the set of vertices of $\bbT_d$. We
single out one vertex from ${\mathcal V}(\bbT_d)$ and call this ${\mathcal O}$, the origin. For each
two vertices $u,v \in {\mathcal V}(\bbT_d)$, we say that $ u \leq v$ if $u$ belongs to the
path connecting ${\mathcal O}$ to $v$. Besides, for two vertices $u,v$ such that $u \leq v$
let $d(u,v)$, be the distance between $u$ and $v$, as the number of edges
the path from $u$ to $v$ has. Now, let us define
\[ \bbT^+_d(u) = \{v \in {\mathcal V}(\bbT_d): u \leq v \}.\]
\noindent
Pick a $v \in {\mathcal V}(\bbT_d)$ such that $d({\mathcal O},v)=1$ and consider 
$\bbT^+_d = \bbT_d \backslash \bbT_d^+(v).$ For $\bbT^+_d$ we define 
\[ \partial T^+_d(u,n) = \{ v \in \bbT^+_d : d(u,v) = n\}.\]
Now we define the \cpp\ in $\bbT_d$. Let $\{R_v\}_{\{ v \in {\mathcal V}(\bbT_d) \}}$ and $R$ 
be a set of independent and identically distributed random variables. We define
$p_k = \bbPr[R = k]$ for $k=0,1,\dots$ To avoid trivialities we assume $p_0 \in (0,1).$
Besides, for each $u \in {\mathcal V}(\bbT_d)$, we define the random sets
\begin{equation}
\label{E: defBu}
B_u = \{v \in {\mathcal V}(\bbT_d): d(u,v) \leq R_u\}
\end{equation}
\noindent
and consider the non-decreasing sequence of random sets $I_0 \subset I_1 \subset \cdots$
defined as $ I_0 = \{{\mathcal O}\} $ and inductively $I_{n+1} = \bigcup_{u \in I_n} B_u$ 
for all $ n \geq 0.$ Let $ I = \bigcup_{n \geq 0} I_n$ be the connected component of the
origin. Under the rumour process interpretation, $I$ is the set of vertices
which heard about the rumour. We say that the process survives if $|I|=\infty,$ 
referring to the surviving event as $V.$

Consider $\bbP_+$ and $\bbP$ the probability measures associated to the processes on $\bbT_d^+$ and $\bbT_d$
(we do not mention the random variable $R$ unless absolutely necessary). By
a coupling argument one can see that for a fixed $R$

\begin{equation}
\bbP_+[V] \leq \bbP[V]
\end{equation}

By the other side, by the definition of $\bbT_d^+$ and its relation with $\bbT_d$ we have that
for a fixed $R$

\begin{equation}
\label{E: Equal}
\bbP_+[V]=0 \hbox{ if and only if } \bbP[V]=0.
\end{equation}

The paper is organized as follows. Section~\ref{S: MR} presents the main results.
Section~\ref{S: proofs} brings the proofs for the main results together with
auxiliary lemmas and handy inequalities. Section~\ref{S: Heterogeneous} presents
results for the heterogeneous setup of the \cpp. Finally, in Section~\ref{S: Exa} 
we present examples where some conditions can be verified.

\section{Main Results}
\label{S: MR}

\begin{teo} 
\label{T:CSPAH}
Consider the \cpp\ on $\bbT_d^+$ with radius of influence $R$
\begin{enumerate}
\item[\textit{(I)}] If $\bbE(d^R) > 1 + p_0$ then, $\bbP_+[V] > 0,$
\item[\textit{(II)}] If $\bbE(d^R) \leq 2 - \frac{1}{d}$ then, $\bbP_+[V] = 0.$
\end{enumerate}
\end{teo}

Let $\rho$ and $\psi$ be, respectively, the smallest non-negative root of the equations
\begin{align}
& \bbE(\rho^{d^R}) + (1 - \rho)p_0 = \rho, \\
& \bbE(\psi^{\frac{d}{d-1}(d^{R}-1)}) = \psi.
\end{align}

\begin{teo} 
\label{T: SobrevivenciaTd+}
Consider the \cpp\ on $\bbT_d^+.$
Then,
\begin{displaymath}
1 - \rho \leq \bbP_+(V) \leq 1 - \psi.
\end{displaymath}
\end{teo}

\begin{teo} 
\label{T: ViveTd}
For the \cpp\ on $\bbT_d$ with radius of influence $R$, it holds that
\begin{equation}
\label{E: ViveTd}
1 - \displaystyle \left(1 -
\rho^{\frac{d+1}{d}}\right)p_0 -
\bbE\displaystyle \left(\rho^{\frac{(d+1)}{d}d^{R}}\right)
\leq \bbP[V] \leq 1 - \bbE\displaystyle
\left(\psi^{\frac{(d+1)}{d-1}(d^{R}-1)}\right).
\end{equation}
\end{teo}

\section{Proofs}
\label{S: proofs}

\subsection{Auxiliary Processes}
\label{SS: AP}

Let us define two auxiliary branching process, being the first one 
$\{\mathcal{X}_n\}_{n \in \bbN}$. For this process, the associated random
variable is  $ X,$ assuming values in $\{0, d, d^2, \dots \} $ 
such that 
\begin{align*}
&\bbPr[X=0] = p_o, \\
&\bbPr[X=d^k] = p_k \hbox{ for } k =1,2, \dots
\end{align*}

\noindent
whose expectation is 
\begin{align}
\label{E: Minorante}
\bbE[X] = \bbE[d^R]-p_0
\end{align}

\noindent
and whose generating function is
\begin{align}
\label{E: GFMinorante}
\varphi_{X}(s) = \bbE[s^X] = \bbE[s^{d^R}] + (1-s)p_0.
\end{align}

The second auxiliary process is $\{\mathcal{Y}_n\}_{n \in \bbN}$. 
For this process, the associated random
variable is  $ Y,$ assuming values in $\{0, d, d + d^2, \dots, \sum_{i=1}^k d^i \} $ 
such that 
\begin{align*}
&\bbPr\Big[Y = \frac{d(d^k-1)}{d-1}\Big] = p_k \hbox{ for } k = 0,1,2, \dots
\end{align*}

\noindent
whose expectation is 
\begin{align}
\label{E: Majorante}
\bbE[Y] = \frac{d}{d-1}(\bbE[d^R]-1)
\end{align}

\noindent
and whose generating function is
\begin{align}
\label{E: GFMajorante}
\varphi_{Y}(s) = \bbE[s^Y] = \bbE[s^{\frac{d}{d-1}(d^{R}-1)}].
\end{align}

\subsection{Proofs}
\label{SS: Proofs}
\hfill

\noindent 
\textit{Proof of Theorem~\ref{T:CSPAH}}\\
By a coupling argument one can see that our process dominates
$\{\mathcal{X}_n\}_{n \in \bbN}$. This process survives as long as
$\bbE[X]> 1$ therefore from~(\ref{E: Minorante}) our process survives
if $\bbE[d^R] > 1 + p_0,$ proving \textit{(I)}. 

By the other side, also by a coupling argument, our process is do\-mi\-na\-ted by
$\{\mathcal{Y}_n\}_{n \in \bbN}.$ That process dies out provided $\bbE[Y] \leq 1$
therefore from~(\ref{E: Majorante}) our process dies out if 
$\bbE[d^R] \leq 2 - \frac{1}{d},$ proving \textit{(II)}.\qed

\noindent 
\textit{Proof of Theorem~\ref{T: SobrevivenciaTd+}}

In order to find the extinction probability of $\{\mathcal{X}_n\}_{n \in \bbN}$ 
(Grimmett and Stirzaker(~\cite[p.173]{GrimmettStirzaker}),
let us consider the smallest non-negative root of the equation $\rho = \varphi_{X}(\rho).$
Therefore from~(\ref{E: GFMinorante})
\[ \bbE[\rho^{d^R}] + (1-\rho)p_0 = \rho \]
and by construction of the processes, as $ \bbP_+[V^c] \leq \rho, $ we have that 
\[ 1- \rho \leq \bbP_+[V] .\]
In order to find the extinction probability of $\{\mathcal{Y}_n\}_{n \in \bbN}$ 
(Grimmett and Stirzaker(~\cite[p.173]{GrimmettStirzaker}), 
let us consider the smallest non-negative root of the equation $\psi = \varphi_{Y}(\psi).$
Therefore from~(\ref{E: GFMajorante})
\[ \bbE[\psi^{\frac{d}{d-1}(d^{R}-1)}] = \psi \]
and by the construction of the processes, as $ \bbP_+[V^c] \geq \psi$, we have that 
\[ \bbP_+[V] \leq 1 - \psi.\]
\qed

\noindent 
\textit{Proof of Theorem~\ref{T: ViveTd}}\\
Observe that except for the root, all vertices see towards 
infinity a tree like $\bbT_d^+.$ So, assuming $R_0=k$ the probability
for the process to survive is larger or equal than the probability of the
process to survive from at least one of the $d^{k-1}(d+1)$ trees that have
as root the furthest infected vertices. By the other side, still assuming
$R_0=k$, the probability for the process to survive in $\bbT_d$  is smaller
or equal than the probability for the process to survive from at least one 
of the $(d+1)(d^k-1)(d-1)^{-1}$ trees like $\bbT_d^+$ that are seen from 
each active vertices by its own, independently. So, for $k = 1,2,\dots$ 
\[
1 - (1 - \bbP_+[V])^{(d+1)d^{k-1}} \leq \bbP[V | R_0 =k] \leq 
1 - (1 - \bbP_+[V])^{\frac{(d+1)}{d-1}[d^{k}-1]}.
\]
From this and from Theorem~\ref{T: SobrevivenciaTd+} follows~(\ref{E: ViveTd}).
\qed

\section{Heterogeneous Cone Percolation on $\bbT_d^+$}
\label{S: Heterogeneous}

Suppose we have two sets of independent random variables, $\{R_z\}_{\{ z \in \bbN \}}$ and
$\{{\bar R}_v\}_{\{ v \in {\mathcal V}(\bbT^+_d) \}},$ such that for all $z \in \bbN$ and all 
$u \in {\mathcal{V}}$ such that $d({\mathcal O},u) = z,$ ${\bar R}_u$ and $R_z$ are equally 
distributed. Besides assume $\bbP[R_z=0] < 1 $ for all $z \in \bbN.$ 

We define
the \hcpp\ from the set of $B_u$ defined in~(\ref{E: defBu}).
For $n \in \bbN$ fixed and $u \leq v \in {\mathcal V}(\bbT^+_d)$, consider the event
\[ V^n_{u,v}: \textit{Process starting from $u$ reaches $v$ in at most $n$ steps}.\]

For a fixed integer $n,$ let $ X_0^n = \{{\mathcal O}\}.$ Besides, for 
$ j=1,2, \dots $ we define
\[ X_j^n = \bigcup_{u \in X_{j-1}^n} \{ v \in \partial
T_n^u : V^n_{u,v} \textit{ occurs } \}.\]
Again, for all $ j=1,2, \dots $ consider
\[Z_j^n = | X_j^n |.\]

So, for all fixed positive integer $n$, $\{Z_j^n\}_{j \geq 0}$ is a branching process dominated
by the number of vertices $v \in \partial T^{\mathcal O}_{jn}$ which are activated.
\begin{lem} 
\label{aes1}
Consider $n$ fixed. For $\mu_j := \bbE[Z^n_j],$ the mean number of offspring on generation $j$ for the process 
$\{Z_j^n\}_{j \geq 0}$, it holds that
\begin{displaymath}
\mu_j = d^n \rho_{j}^n,
\end{displaymath}
where $\rho_{j}^n = \bbP[V^n_{u,v}]$, for any fixed pair $u \leq v$ such that $d({\mathcal O},u) =jn$ 
and $d({\mathcal O},v)=(j+1)n.$
\end{lem}
\noindent \textit{Proof of Lemma~\ref{aes1}}\\
For fixed $j$ and $n$, consider $\partial T_n^v = \{ u_1, u_2, \dots , u_{d^n} \}$. So we can write
$Z_j^n$ as $\sum_{i=1}^{d^n}I_{\{V^n_{v,u_i} \}}.$ Taking expectation in both sides finishes the proof.

\begin{lem} \label{aes2}
Consider $n$ fixed and $\rho_{j}^n = \bbP[V^n_{u,v}]$, for any fixed pair $u \leq v$ such that 
$d({\mathcal O},u) =jn$ and $d({\mathcal O},v)=(j+1)n.$
\begin{displaymath}
\rho_{j}^n \geq
\prod_{k=0}^{n-1}[1-\prod_{i=0}^{k}\bbP[R_{jk+i} < k+1-i]].
\end{displaymath}
\end{lem}
\noindent \textit{Proof of Lemma~\ref{aes2}}\\
For any fixed pair $u \leq v$ such that $d({\mathcal O},u) =jn$ and $d({\mathcal O},v)
=(j+1)n$ we have that

\begin{displaymath}
V^n_{u,v} = \bigcap_{k=0}^{n-1}\displaystyle
\left[\bigcup_{i=0}^{k}\{ R_{jn+i} \geq k+1-i \}\right]
\end{displaymath}
\noindent
and so

\begin{align*}
\rho_{j}^n = & \ \bbP\displaystyle \left (
\bigcap_{k=0}^{n-1}\displaystyle \left[\bigcup_{i=0}^{k}\{
R_{jn+i} \geq k+1-i \}\right] \right) \\ \geq & \ 
\prod_{k=0}^{n-1}\bbP\displaystyle \left(\bigcup_{i=0}^{k}\{
R_{jn+i} \geq k+1-i \}\right).
\end{align*}
The inequality is a consequence of the FKG inequality (Grimmett~\cite[p.34]{Grimmett}). \qed

\begin{teo} \label{aest}
The \hcpp\ in $\bbT^+_d$ has a giant component with positive pro\-ba\-bi\-li\-ty if
for some fixed $n$,
\begin{equation}
\label{eq1}
\liminf_{j \rightarrow \infty} d^n
\prod_{k=0}^{n-1}[1-\prod_{i=0}^{k}\bbP[R_{jk+i} < k+1-i]] >
1.
\end{equation}
\end{teo}
\noindent \textit{Proof of Theorem~\ref{aest}}\\
From Theorem~1 in Souza \& Biggins~(\cite[p.39]{SouzaBiggins}), a branching process in varying environments
is uniformly supercritical if
\begin{displaymath}
\liminf_{j \rightarrow \infty} \mu_j > 1.
\end{displaymath}
From Lemma~\ref{aes1} and Lemma~\ref{aes2}, that is what happens if~(\ref{eq1}) holds.

From the fact that
\[\frac{Z_j^n}{\bbE[Z_j^n]} \leq
\frac{1}{\rho_{j}^n} \]
\noindent
one can see that the \hcpp\ has a giant component with positive probability if
\begin{displaymath}
\liminf_{j \rightarrow \infty} d^n
\prod_{k=0}^{n-1}[1-\prod_{i=0}^{k}\bbP[R_{jk+i} < k+1-i]] >
1.
\end{displaymath}
\qed

\section{Examples}
\label{S: Exa}

\begin{exa}
\label{E: B(p)}
Consider a \cpp\ in $\bbT_d,$ assuming
\[
\bbPr[R = 1] = p = 1 - \bbPr[R = 0].
\]
In words $R \sim {\mathcal B}(p).$

\begin{itemize}
\item If $p > d^{-1}$ then, $\bbP[V]>0.$
\item If $p \leq d^{-1}$ then, $\bbP[V]=0.$
\end{itemize}

By the definition, one can see that
\[
\bbP[V^c] = (1-p) + {(\bbP_+[V^c])}^{d+1}p.
\]
\end{exa}

Observing that the upper and lower process presented by
$\{\mathcal{X}_n\}_{n \in \bbN}$ and $\{\mathcal{Y}_n\}_{n \in \bbN}$ 
presented in session~\ref{SS: AP} are the same, we see that
\[
\bbP[V]=p(1-\psi^{d+1}),
\]
being $\psi$ the solution of
\[
p\psi^d-\psi+1-p=0.
\]

\begin{exa} 
\label{E: G1-p}
Consider a \cpp\ in $\bbT_d,$ assuming
\[
\bbPr(R = k) = (1-p)p^{k}, k=0,1,2,\dots
\]

In other words $R \sim {\mathcal G}(1-p).$ From Theorem~\ref{T: ViveTd}

\begin{itemize}
\item If $dp^2 - 2dp + 1 < 0$ then, $\bbP[V]>0$.
\item If $2pd \leq 1$ then, $\bbP[V]=0$.
\end{itemize}

As a consequence of this and~(\ref{E: Equal}), for $d$ fixed
\[
\frac{1}{2d} < \inf\{p: \bbP[V]>0\} \leq 1 - \sqrt{1- \frac{1}{d}}.
\]
\end{exa}

\begin{exa} 
\label{C:CSPAHB}
Consider a \cpp\ in $\bbT_d,$ assuming
\[
\bbPr(R = k) = \binom{n}{k}p^k(1-p)^{n-k}  , \quad k=0,1,\dots, n
\]

\begin{itemize}
\item If $(pd+1-p)^n - (1-p)^n > 1$ then, $\bbP[V]>0$.
\item If $2d - d(pd+1-p)^n \geq 1$ then, $\bbP[V]=0$.
\end{itemize}

Let $d=2$ and $R \sim \mathcal{B}(4,\frac{1}{2})$. 

Therefore $\rho$ and $\psi$ are, respectively,
solutions of
\begin{align*}
& x^{16} + 4x^8 + 6x^4 + 4x^2 - 16x + 1 = 0, \\& x^{30} + 4x^{14} +
6x^6 + 4x^2 - 16x + 1 = 0.
\end{align*}
So $\rho$ = 0.0635146 and $\psi$ = 0.06350850, which implies that
\[
0.937435919 \leq \bbP[V] \leq 0.937435962.
\]

Let $d=4$ and $R \sim \mathcal{B}(4,\frac{1}{4})$. 

Therefore $\rho$ and $\psi$ are, respectively,
solutions of
\begin{align*}
& x^{256} + 12x^{64} + 54x^{16} + 108x^{4} - 256x + 81 = 0, \\&
x^{340} + 12x^{84} + 54x^{20} + 108x^{4} - 256x + 81 = 0.
\end{align*}
So $\rho$ = 0.3208787235 and $\psi$ = 0.3208787200, which implies that
\[
0.682158629 \leq \bbP[V] \leq 0,682158630.
\]
\end{exa}

\begin{exa} \label{aeslp}
Consider a \hcpp\ on $\bbT^+_d,$ assuming that $R_j$ are Bernoullis, that is,
\begin{displaymath}
\bbP[R_j = 1] = 1 - \bbP[R_j = 0] \hbox{ for } j = 0,1,2\dots
\end{displaymath}
By applying Theorem~\ref{aest} with $n = 1$ one can see that the \hcpp\ 
on $\bbT^+_d$ survives with positive probability if
\begin{displaymath}
\liminf_{j\rightarrow \infty} d \bbP[R_j = 1] > 1.
\end{displaymath}

\end{exa}

\end{document}